\documentclass[10pt]{article}
\usepackage{cite}
\usepackage{mathrsfs}
\usepackage{amsfonts}
\usepackage{amsmath}
\usepackage{amsfonts,amssymb,color}
\usepackage{dsfont}
\usepackage{curves}
\usepackage{mathrsfs}
\usepackage{pifont}
\usepackage{amssymb}
\allowdisplaybreaks

\numberwithin{equation}{section}

\date{}

\def\BigRoman{\uppercase\expandafter{\romannumeral\number\count 255 }}
\def\Romannumeral{\afterassignment\BigRoman\count255=}

\begin{document}
\title{Spanning $k$-trees, odd $[1,b]$-factors and spectral radius in binding graphs
}
\author{\small  Jiancheng Wu, Sizhong Zhou\footnote{Corresponding author. E-mail address: zsz\_cumt@163.com (S. Zhou)}\\
\small  School of Science, Jiangsu University of Science and Technology,\\
\small  Zhenjiang, Jiangsu 212100, China\\
}

\maketitle
\begin{abstract}
\noindent The binding number of a graph $G$, written as $\mbox{bind}(G)$, is defined by
$$
\mbox{bind}(G)=\min\left\{\frac{|N_G(X)|}{|X|}:\emptyset\neq X\subseteq V(G),N_G(X)\neq V(G)\right\}.
$$
A graph $G$ is called $r$-binding if $\mbox{bind}(G)\geq r$. An odd $[1,b]$-factor of a graph $G$ is a spanning subgraph $F$ with $d_F(v)\in\{1,3,\ldots,b\}$ for all $v\in V(G)$, where $b\geq1$
is an odd integer. A spanning $k$-tree of a connected graph $G$ is a spanning tree $T$ with $d_T(v)\leq k$ for every $v\in V(G)$. In this paper, we first show a tight sufficient condition with
respect to the adjacency spectral radius for connected $\frac{1}{b}$-binding graphs to have odd $[1,b]$-factors, which generalizes Fan and Lin's previous result [D. Fan, H. Lin, Binding number,
$k$-factor and spectral radius of graphs, Electron. J. Combin. 31(1) (2024) \#P1.30] and partly improves Fan, Liu and Ao's previous result [A. Fan, R. Liu, G. Ao, Spectral radius, odd $[1,b]$-factor
and spanning $k$-tree of 1-binding graphs, Linear Algebra Appl. 705 (2025) 1--16]. Then we put forward a tight sufficient condition via the adjacency spectral radius for connected $\frac{1}{k-2}$-binding
graphs to have spanning $k$-trees, which partly improves Fan, Liu and Ao's previous result [A. Fan, R. Liu, G. Ao, Spectral radius, odd $[1,b]$-factor and spanning $k$-tree of 1-binding graphs,
Linear Algebra Appl. 705 (2025) 1--16].
\\
\begin{flushleft}

(2020) Mathematics Subject Classification: 05C50, 05C05, 05C70
\end{flushleft}
\end{abstract}

\section{Introduction}

In this paper, we only deal with undirected and simple graphs. Let $G$ denote a graph with vertex set $V(G)$ and edge set $E(G)$, where $e(G)=|E(G)|$ and $n=|V(G)|$ are called the size and
the order of $G$, respectively. We denote by $c(G)$ the number of components of $G$, and by $o(G)$ the number of odd components of $G$. For $v\in V(G)$, the set of vertices adjacent to $v$
in $G$, denoted by $N_G(v)$, is called the neighborhood of $v$ in $G$. For $v\in V(G)$, $d_G(v)=|N_G(v)|$ is called the degree of $v$ in $G$. For any $S\subseteq V(G)$, the subgraph of $G$
induced by $S$ is denoted by $G[S]$, and write $G-S=G[V(G)\setminus S]$. Let $G_1$ and $G_2$ denote two vertex-disjoint graphs. The union of $G_1$ and $G_2$ is denoted by $G_1\cup G_2$. The
join $G_1\vee G_2$ is obtained from $G_1\cup G_2$ by adding all possible edges between $V(G_1)$ and $V(G_2)$. For any integer $t\geq1$, we denote by $tG$ the union of $t$ copies of $G$.

Woodall \cite{Wt} introduced the concept of the binding number. For any subset $X\subseteq V(G)$, we write $N_G(X)=\bigcup\limits_{x\in X}N_G(x)$. The binding number of a graph $G$, written
as $\mbox{bind}(G)$, is defined by
$$
\mbox{bind}(G)=\min\left\{\frac{|N_G(X)|}{|X|}:\emptyset\neq X\subseteq V(G),N_G(X)\neq V(G)\right\}.
$$
A graph $G$ is called $r$-binding if $\mbox{bind}(G)\geq r$.

Given a graph $G$ with $V(G)=\{v_1,v_2,\ldots,v_n\}$, the adjacency matrix $A(G)=(a_{ij})$ of $G$ is an $n\times n$ symmetric matrix, where $a_{ij}=1$ if $v_iv_j\in E(G)$ and $a_{ij}=0$ if
$v_iv_j\notin E(G)$. The largest eigenvalue of $A(G)$, denoted by $\rho(G)$, is referred to as the adjacency spectral radius of $G$.

Let $a$ and $b$ be two integers with $1\leq a\leq b$. A spanning subgraph $F$ of $G$ satisfying $d_F(v)\in[a,b]$ for any $v\in V(G)$ is called an $[a,b]$-factor of $G$. An odd $[1,b]$-factor
of $G$ is a spanning subgraph $F$ with $d_F(v)\in\{1,3,\ldots,b\}$ for all $v\in V(G)$, where $b\geq1$ is an odd integer. In particular, an odd $[1,b]$-factor is a 1-factor (or a perfect
matching) if $b=1$.

Lots of scholars investigated the existence of factors in graphs from the prospective of binding number, see \cite{Ap,RW,KT,GW,KW}. Many scholars focused on the existence of factors in graphs
from the prospective of the spectral radius, see \cite{O,LM,Zs,ZW1,ZZS,ZSZ,Zt,ZZ,ZZL,ZWH,FLL,ZL,WZ}. Fan and Lin \cite{FL} provided an adjacency spectral condition for a connected 1-binding graph
to contain a 1-factor. Fan, Liu and Ao \cite{FLA} arose a tight sufficient condition with respect to the adjacency spectral radius to guarantee that a connected 1-binding graph has an odd
$[1,b]$-factor.

\medskip

\noindent{\textbf{Theorem 1.1}} (Fan, Liu and Ao \cite{FLA}). Let $G$ be a connected 1-binding graph of even order $n\geq18b+6$, where $b\geq1$ is an odd integer. If
$$
\rho(G)\geq\rho(K_1\vee(K_{n-3b-2}\cup bK_3\cup K_1)),
$$
then $G$ has an odd $[1,b]$-factor, unless $G=K_1\vee(K_{n-3b-2}\cup bK_3\cup K_1)$.

\medskip

Motivated by \cite{A,FLA} directly, we also establish a connection between the adjacency spectral radius and odd $[1,b]$-factors in binding graphs, and obtain the following result.

\medskip

\noindent{\textbf{Theorem 1.2.}} Let $b\geq1$ be an odd integer, and let $G$ be a connected $\frac{1}{b}$-binding graph of even order $n\geq\max\{12,2b+8\}$. If
$$
\rho(G)\geq\rho(K_1\vee(K_{n-b-4}\cup K_3\cup bK_1)),
$$
then $G$ has an odd $[1,b]$-factor, unless $G=K_1\vee(K_{n-b-4}\cup K_3\cup bK_1)$.

\medskip

One easily checks that $\rho(K_1\vee(K_{n-3b-2}\cup bK_3\cup K_1))\leq\rho(K_1\vee(K_{n-b-4}\cup K_3\cup bK_1))$. Hence, the spectral condition in Theorem 1.1 is better than that of Theorem 1.2.
But under the spectral condition $\rho(G)\geq\rho(K_1\vee(K_{n-b-4}\cup K_3\cup bK_1))$, our order and binding number conditions are better than ones of Theorem 1.1. Furthermore, our improvements
on order and binding number in Theorem 1.2 are very large. From the above discussion, Theorem 1.1 and Theorem 1.2 cannot be obtained from each other, and so this large improvements are interesting
and important.

A spanning $k$-tree of a connected graph $G$ is a spanning tree $T$ with $d_T(v)\leq k$ for every $v\in V(G)$. Lots of scholars paid much attention to the existence of spanning trees in connected
graphs from the prospective of the spectral radius, see \cite{ALY,Wc,ZSL,ZZL1,ZW}. Fan, Goryainov, Huang and Lin \cite{FGHL} showed a spectral characterization for a connected graph with a spanning
$k$-tree, where $k\geq3$ is an integer. Fan, Liu and Ao \cite{FLA} provided a tight sufficient condition with respect to the adjacency spectral radius to guarantee that a connected 1-binding graph
has a spanning $k$-tree, where $k\geq8$ is an integer.

\medskip

\noindent{\textbf{Theorem 1.3}} (Fan, Liu and Ao \cite{FLA}). Let $k\geq8$, and let $G$ be a connected 1-binding graph of order $n\geq k^{3}-4k^{2}+5k+4$. If
$$
\rho(G)\geq\rho(K_1\vee(K_{n-2k}\cup(k-1)K_2\cup K_1)),
$$
then $G$ has a spanning $k$-tree, unless $G=K_1\vee(K_{n-2k}\cup(k-1)K_2\cup K_1)$.

\medskip

Motivated by \cite{W,FLA} directly, we also establish a relationship between the adjacency spectral radius and spanning $k$-trees in binding graphs, and put forward the following result.

\medskip

\noindent{\textbf{Theorem 1.4.}} Let $k\geq3$ be an integer, and let $G$ be a connected $\frac{1}{k-2}$-binding graph of order $n\geq2k+22$. If
$$
\rho(G)\geq\rho(K_1\vee(K_{n-k-3}\cup2K_2\cup(k-2)K_1)),
$$
then $G$ has a spanning $k$-tree, unless $G=K_1\vee(K_{n-k-3}\cup2K_2\cup(k-2)K_1)$.

\medskip

Obviously, the range of the value of the parameter $k$ in Theorem 1.4 is wider than that in Theorem 1.3. One easily checks that $\rho(K_1\vee(K_{n-2k}\cup(k-1)K_2\cup K_1))\leq\rho(K_1\vee(K_{n-k-3}\cup2K_2\cup(k-2)K_1))$.
Consequently, the spectral condition in Theorem 1.3 is better than that of Theorem 1.4. But under the spectral condition $\rho(G)\geq\rho(K_1\vee(K_{n-k-3}\cup2K_2\cup(k-2)K_1))$, our order and binding
number conditions are better than ones of Theorem 1.3. Furthermore, our improvements on parameter $k$, order and binding number in Theorem 1.4 are very large. From the above discussion, Theorem 1.3 and
Theorem 1.4 cannot be obtained from each other, and so this large improvements are interesting and important.

\section{Preliminary lemmas}

In this section, we provide some useful lemmas, which will be used to prove our main results. Amahashi \cite{A} showed a characterization for a graph with an odd $[1,b]$-factor, where $b\geq1$ is
an odd integer.

\medskip

\noindent{\textbf{Lemma 2.1}} (Amahashi \cite{A}). Let $b\geq1$ be an odd integer and let $G$ be a graph. Then $G$ has an odd $[1,b]$-factor if and only if
$$
o(G-S)\leq b|S|
$$
for all $S\subseteq V(G)$.

\medskip

Win \cite{W} provided a sufficient condition for a connected graph to possess a spanning $k$-tree.

\medskip

\noindent{\textbf{Lemma 2.2}} (Win \cite{W}). Let $k\geq2$ be an integer and let $G$ be a connected graph. Then $G$ has a spanning $k$-tree if
$$
c(G-S)\leq(k-2)|S|+2
$$
for any $S\subseteq V(G)$.

\medskip

\noindent{\textbf{Lemma 2.3}} (Li and Feng \cite{LF}). Let $G$ be a connected graph and $H\subseteq G$. Then
$$
\rho(G)\geq\rho(H),
$$
with equality occurring if and only if $G=H$.

\medskip

\noindent{\textbf{Lemma 2.4}} (Hong \cite{Ha}). Let $G$ be a graph of order $n$. Then
$$
\rho(G)\leq\sqrt{2e(G)-n+1},
$$
where the equality occurs if and only if $G$ is a star or a complete graph.

\medskip

\noindent{\textbf{Lemma 2.5}} (Fan and Lin \cite{FL}). Let $\sum\limits_{i=1}^{t}n_i=n-s$ with $s\geq1$. If $n_t\geq n_{t-1}\geq\cdots\geq n_1\geq1$, $n_{t-1}\geq3$ and $n_t<n-s-t-1$, then
$$
\rho(K_s\vee(K_{n_t}\cup K_{n_{t-1}}\cup\cdots\cup K_{n_1}))<\rho(K_s\vee(K_{n-s-t-1}\cup K_3\cup(t-2)K_1)).
$$

\medskip

\noindent{\textbf{Lemma 2.6.}} Let $\sum\limits_{i=1}^{t}n_i=n-s$ with $s\geq1$. If $n_t\geq n_{t-1}\geq n_{t-2}\geq\cdots\geq n_1\geq1$, $n_{t-2}\geq2$ and $n_t<n-s-t-1$, then
$$
\rho(K_s\vee(K_{n_t}\cup K_{n_{t-1}}\cup K_{n_{t-2}}\cup\cdots\cup K_{n_1}))<\rho(K_s\vee(K_{n-s-t-1}\cup2K_2\cup(t-3)K_1)).
$$

\medskip

\noindent{\it Proof.} Let $G=K_s\vee(K_{n_t}\cup K_{n_{t-1}}\cup K_{n_{t-2}}\cup\cdots\cup K_{n_1})$, and let $\mathbf{x}$ be the Perron vector of $A(G)$ with respect to $\rho(G)$. According to symmetry, we
let $\mathbf{x}(v)=x_i$ for every $v\in V(K_{n_i})$, where $1\leq i\leq t$, and $\mathbf{x}(u)=y$ for every $u\in V(K_s)$. We easily see that $K_{n_t+s}$ is a proper subgraph of $G$. Combining this with Lemma
2.3, we conclude
$$
\rho(G)>\rho(K_{n_t+s})=n_t+s-1>n_t-1.
$$
Note that $n_t\geq n_i$ for $1\leq i\leq t-1$. In terms of $A(G)\mathbf{x}=\rho(G)\mathbf{x}$, we get
\begin{align*}
(\rho(G)-(n_i-1))(x_t-x_i)=(n_t-n_i)x_t\geq0
\end{align*}
for $1\leq i\leq t-1$. This yields that $x_t\geq x_i$ for $1\leq i\leq t-1$. Let $G'=K_s\vee(K_{n-s-t-1}\cup2K_2\cup(t-3)K_1)$. Then using $2\leq n_{t-2}\leq n_{t-1}\leq n_t<n-s-t-1$, $n_j\geq1$ for
$1\leq j\leq t-3$ and $x_t\geq x_i$ for $1\leq i\leq t-1$, we obtain
\begin{align*}
\rho(G')-\rho(G)\geq&\mathbf{x}^{T}(A(G')-A(G))\mathbf{x}\\
=&2(n_{t-1}-2)x_{t-1}(n_tx_t-2x_{t-1}+(n_{t-2}-2)x_{t-2}+\sum\limits_{j=1}^{t-3}(n_j-1)x_j)\\
&+2(n_{t-2}-2)x_{t-2}(n_tx_t-2x_{t-2})\\
&+2\sum\limits_{j=1}^{t-3}(n_j-1)x_j(n_tx_t-x_j+(n_{t-2}-2)x_{t-2})\\
&+2\sum\limits_{i=2}^{t-3}\sum\limits_{j=i-1}^{t-4}(n_i-1)(n_j-1)x_ix_j\\
>&0,
\end{align*}
which implies that $\rho(G)<\rho(G')$. That is to say, $\rho(K_s\vee(K_{n_t}\cup K_{n_{t-1}}\cup K_{n_{t-2}}\cup\cdots\cup K_{n_1}))<\rho(K_s\vee(K_{n-s-t-1}\cup2K_2\cup(t-3)K_1))$. \hfill $\Box$

\medskip

Let $M$ denote a real $n\times n$ matrix and $\mathcal{N}=\{1,2,\ldots,n\}$. Given a partition $\pi:\mathcal{N}=\mathcal{N}_1\cup\mathcal{N}_2\cup\cdots\cup\mathcal{N}_r$, the matrix $M$ can
be correspondingly given by
\begin{align*}
M=\left(
  \begin{array}{cccc}
    M_{11} & M_{12} & \cdots & M_{1r}\\
    M_{21} & M_{22} & \cdots & M_{2r}\\
    \vdots & \vdots & \ddots & \vdots\\
    M_{r1} & M_{r2} & \cdots & M_{rr}\\
  \end{array}
\right),
\end{align*}
where the block $M_{ij}$ denotes the $n_i\times n_j$ matrix for any $1\leq i,j\leq r$. Let $m_{ij}$ denote the average value of all sums of $M_{ij}$. Then the matrix $M_{\pi}=(m_{ij})_{r\times r}$ is called
the quotient matrix of $M$. The partition $\pi$ is called equitable if every block $M_{ij}$ of $M$ has constant row sum $m_{ij}$ for $1\leq i,j\leq r$.

\medskip

\noindent{\textbf{Lemma 2.7}} (Brouwer and Haemers \cite{BH}, You, Yang, So and Xi \cite{YYSX}). Let $M$ denote a real $n\times n$ matrix with an equitable partition $\pi$, and let $M_{\pi}$ be the corresponding
quotient matrix. Then the eigenvalues of $M_{\pi}$ are eigenvalues of $M$. Furthermore, if $M$ is nonnegative and irreducible, then the largest eigenvalues of $M$ and $M_{\pi}$ are equal.

\medskip

\noindent{\textbf{Lemma 2.8}} (Fan and Lin \cite{FL}). Let $G$ be a connected 1-binding graph of even order $n\geq12$. If
$$
\rho(G)\geq\rho(K_1\vee(K_{n-5}\cup K_3\cup K_1)),
$$
then $G$ has a 1-factor, unless $G=K_1\vee(K_{n-5}\cup K_3\cup K_1)$.

\medskip

\section{The proof of Theorem 1.3}

\noindent{\it Proof of Theorem 1.3.} According to Lemma 2.8, Theorem 1.3 is true for $b=1$. Next, we consider $b\geq3$.

Suppose that a connected $\frac{1}{b}$-binding graph $G$ has no odd $[1,b]$-factor. From Lemma 2.1, there exists a nonempty subset $S$ of $V(G)$ such that $o(G-S)>b|S|$. Since $n$ is even, we conclude
$o(G-S)\equiv b|S|$ (mod 2). Hence, $o(G-S)\geq b|S|+2$. Write $|S|=s$ and $o(G-S)=q$. Then $q\geq bs+2$. We denote by $O_1,O_2,\ldots,O_q$ the $q$ odd components in $G-S$ and write $|O_i|=n_i$ for
$1\leq i\leq q$. Without loss of generality, let $n_q\geq n_{q-1}\geq\cdots\geq n_1$. We are to show $n_{bs+1}\geq3$. Otherwise, we possess $n_i=1$ for $1\leq i\leq bs+1$. Write
$X=V(O_1\cup O_2\cup\cdots\cup O_{bs+1})$. Clearly, $N_G(X)\subseteq S$ and $|X|=bs+1$. Thus, we see
$$
\frac{|N_G(X)|}{|X|}\leq\frac{|S|}{|X|}=\frac{s}{bs+1}<\frac{1}{b},
$$
which is impossible because $G$ is $\frac{1}{b}$-binding. This implies $n_i\geq3$ for $i\geq bs+1$. It is obvious that $G$ is a spanning subgraph of $G_1=K_s\vee(K_{n_1}\cup K_{n_2}\cup\cdots\cup K_{n_{bs+1}}\cup K_{n_{bs+2}})$,
where $n_{bs+2}\geq n_{bs+1}\geq\cdots\geq n_2\geq n_1\geq1$, $n_{bs+1}\geq3$ and $\sum\limits_{i=1}^{bs+2}n_i=n-s$. Using Lemma 2.3, we conclude
\begin{align}\label{eq:3.1}
\rho(G)\leq\rho(G_1),
\end{align}
where the equality occurs if and only if $G=G_1$. Let $G_2=K_s\vee(K_{n-(b+1)s-3}\cup K_3\cup bsK_1)$, where $n\geq(b+1)s+6$. In terms of Lemma 2.5 (set $t=bs+2$ in Lemma 2.5), we infer
\begin{align}\label{eq:3.2}
\rho(G_1)\leq\rho(G_2),
\end{align}
with equality if and only if $G_1=G_2$. For $s=1$, we possess $G_2=K_1\vee(K_{n-b-4}\cup K_3\cup bK_1)$. By virtue of \eqref{eq:3.1} and \eqref{eq:3.2}, we get
$$
\rho(G)\leq\rho(K_1\vee(K_{n-b-4}\cup K_3\cup bK_1)),
$$
with equality occurring if and only if $G=K_1\vee(K_{n-b-4}\cup K_3\cup bK_1)$. Observe that $K_1\vee(K_{n-b-4}\cup K_3\cup bK_1)$ has no odd $[1,b]$-factor, a contradiction. We shall consider $s\geq2$.

Recall that $G_2=K_s\vee(K_{n-(b+1)s-3}\cup K_3\cup bsK_1)$. By the partition $V(G_2)=V(K_s)\cup V(K_{n-(b+1)s-3})\cup V(K_3)\cup V(bsK_1)$, the quotient matrix of $A(G_2)$ equals
\begin{align*}
B_2=\left(
  \begin{array}{cccc}
  s-1 & n-(b+1)s-3 & 3 & bs\\
  s & n-(b+1)s-4 & 0 & 0\\
  s & 0 & 2 & 0\\
  s & 0 & 0 & 0\\
  \end{array}
\right).
\end{align*}
The characteristic polynomial of $B_2$ is
\begin{align}\label{eq:3.3}
\varphi(B_2,x)=&x^{4}+(-n+bs+3)x^{3}+(n-bs^{2}-(b+3)s-6)x^{2}\nonumber\\
&+(bs^{2}n+3sn+2n-b(b+1)s^{3}-(5b+3)s^{2}-(2b+12)s-8)x\nonumber\\
&-2bs^{2}n+2b(b+1)s^{3}+8bs^{2}.
\end{align}
Clearly, the partition $V(G_2)=V(K_s)\cup V(K_{n-(b+1)s-3})\cup V(K_3)\cup V(bsK_1)$ is equitable. From Lemma 2.7, $\rho(G_2)$ is the largest root of $\varphi(B_2,x)=0$.

Let $G_*=K_1\vee(K_{n-b-4}\cup K_3\cup bK_1)$. Obviously, $G_*=G_2$ when $s=1$. Let $B_*$ denote the quotient matrix of $A(G_*)$. According to \eqref{eq:3.3}, we obtain
\begin{align*}
\varphi(B_*,x)=&x^{4}+(-n+b+3)x^{3}+(n-2b-9)x^{2}\\
&+(bn+5n-b^{2}-8b-23)x-2bn+2b^{2}+10b.
\end{align*}
By Lemma 2.7, $\rho(G_*)$ is the largest root of $\varphi(B_*,x)=0$. Notice that $K_{n-b-3}$ is a proper subgraph of $G_*=K_1\vee(K_{n-b-4}\cup K_3\cup bK_1)$. In view of Lemma 2.3, we get
\begin{align}\label{eq:3.4}
\rho(G_*)>\rho(K_{n-b-3})=n-b-4.
\end{align}
A simple computation yields that
\begin{align}\label{eq:3.5}
\varphi(B_*,x)-\varphi(B_2,x)=&(s-1)(-bx^{3}+(bs+2b+3)x^{2}-(bsn+bn+3n\nonumber\\
&-(b^{2}+b)s^{2}-(b^{2}+6b+3)s-b^{2}-8b-15)x+2bsn+2bn\nonumber\\
&-(2b^{2}+2b)s^{2}-(2b^{2}+10b)s-2b^{2}-10b).
\end{align}
Let $f(x)=-bx^{3}+(bs+2b+3)x^{2}-(bsn+bn+3n-(b^{2}+b)s^{2}-(b^{2}+6b+3)s-b^{2}-8b-15)x+2bsn+2bn-(2b^{2}+2b)s^{2}-(2b^{2}+10b)s-2b^{2}-10b$. Recall that $n\geq(b+1)s+6$. For $x\geq n-b-4$, it
follows from $b\geq3$, $s\geq2$ and $n\geq(b+1)s+6$ that
\begin{align*}
f'(x)=&-3bx^{2}+2(bs+2b+3)x-bsn-bn-3n+(b^{2}+b)s^{2}+(b^{2}+6b+3)s+b^{2}+8b+15\\
\leq&-3b(n-b-4)^{2}+2(bs+2b+3)(n-b-4)-bsn-bn-3n\\
&+(b^{2}+b)s^{2}+(b^{2}+6b+3)s+b^{2}+8b+15\\
=&-3bn^{2}+(bs+6b^{2}+27b+3)n+(b^{2}+b)s^{2}\\
&-(b^{2}+2b-3)s-3b^{3}-27b^{2}-62b-9\\
\leq&-3b((b+1)s+6)^{2}+(bs+6b^{2}+27b+3)((b+1)s+6)\\
&+(b^{2}+b)s^{2}-(b^{2}+2b-3)s-3b^{3}-27b^{2}-62b-9\\
=&-(3b^{3}+4b^{2}+b)s^{2}+(6b^{3}-4b^{2}-2b+6)s-3b^{3}+9b^{2}-8b+9\\
\leq&-4(3b^{3}+4b^{2}+b)+2(6b^{3}-4b^{2}-2b+6)-3b^{3}+9b^{2}-8b+9\\
=&-3b^{3}-15b^{2}-16b+21\\
<&0,
\end{align*}
which yields that $f(x)$ is decreasing with respect to $x\geq n-b-4$. Hence, we admit
\begin{align}\label{eq:3.6}
f(x)\leq&f(n-b-4)\nonumber\\
=&-b(n-b-4)^{3}+(bs+2b+3)(n-b-4)^{2}\nonumber\\
&-(bsn+bn+3n-(b^{2}+b)s^{2}-(b^{2}+6b+3)s-b^{2}-8b-15)(n-b-4)\nonumber\\
&+2bsn+2bn-(2b^{2}+2b)s^{2}-(2b^{2}+10b)s-2b^{2}-10b\nonumber\\
=&-bn^{3}+(3b^{2}+13b)n^{2}+((b^{2}+b)s^{2}+(4b+3)s-3b^{3}-26b^{2}-53b+3)n\nonumber\\
&-(b^{3}+7b^{2}+6b)s^{2}-(4b^{2}+21b+12)s+b^{4}+13b^{3}+53b^{2}+63b-12.
\end{align}
Let $g(n)=-bn^{3}+(3b^{2}+13b)n^{2}+((b^{2}+b)s^{2}+(4b+3)s-3b^{3}-26b^{2}-53b+3)n-(b^{3}+7b^{2}+6b)s^{2}-(4b^{2}+21b+12)s+b^{4}+13b^{3}+53b^{2}+63b-12$. For $n\geq(b+1)s+6$, it follows from $b\geq3$
and $s\geq2$ that
\begin{align*}
g'(n)=&-3bn^{2}+2(3b^{2}+13b)n+(b^{2}+b)s^{2}+(4b+3)s-3b^{3}-26b^{2}-53b+3\\
\leq&-3b((b+1)s+6)^{2}+2(3b^{2}+13b)((b+1)s+6)\\
&+(b^{2}+b)s^{2}+(4b+3)s-3b^{3}-26b^{2}-53b+3\\
=&-(3b^{3}+5b^{2}+2b)s^{2}+(6b^{3}-4b^{2}-6b+3)s-3b^{3}+10b^{2}-5b+3\\
\leq&-4(3b^{3}+5b^{2}+2b)+2(6b^{3}-4b^{2}-6b+3)-3b^{3}+10b^{2}-5b+3\\
=&-3b^{3}-18b^{2}-25b+9\\
<&0,
\end{align*}
which implies that $g(n)$ is decreasing with respect to $n\geq(b+1)s+6$. For $n\geq(b+1)s+6$, we obtain
\begin{align}\label{eq:3.7}
g(n)\leq&g((b+1)s+6)\nonumber\\
=&-b^{2}(b+1)^{2}s^{3}+(b+1)(3b^{3}-3b^{2}-b+3)s^{2}\nonumber\\
&+(-3b^{4}+7b^{3}+b^{2}+b+9)s+b^{4}-5b^{3}+5b^{2}-3b+6.
\end{align}
Let $h(s)=-b^{2}(b+1)^{2}s^{3}+(b+1)(3b^{3}-3b^{2}-b+3)s^{2}+(-3b^{4}+7b^{3}+b^{2}+b+9)s+b^{4}-5b^{3}+5b^{2}-3b+6$. For $s\geq2$, we deduce
\begin{align*}
h'(s)=&-3b^{2}(b+1)^{2}s^{2}+2(b+1)(3b^{3}-3b^{2}-b+3)s-3b^{4}+7b^{3}+b^{2}+b+9\\
\leq&-12b^{2}(b+1)^{2}+4(b+1)(3b^{3}-3b^{2}-b+3)-3b^{4}+7b^{3}+b^{2}+b+9\\
=&-3b^{4}-17b^{3}-27b^{2}+9b+21\\
<&0 \ \ \ \ \ (\mbox{since} \ b\geq3).
\end{align*}
It follows that $h(s)$ is decreasing with respect to $s\geq2$. For $s\geq2$ and $b\geq3$, we deduce
\begin{align}\label{eq:3.8}
h(s)\leq h(2)=-b^{4}-7b^{3}-17b^{2}+7b+36<0.
\end{align}

By \eqref{eq:3.5}, \eqref{eq:3.6}, \eqref{eq:3.7}, \eqref{eq:3.8} and $s\geq2$, we infer
$$
\varphi(B_*,x)-\varphi(B_2,x)=(s-1)f(x)\leq(s-1)g(n)\leq(s-1)h(s)<0,
$$
which yields
$$
\varphi(B_2,x)>\varphi(B_*,x)
$$
for $x\geq n-b-4$. Combining this with \eqref{eq:3.1} and \eqref{eq:3.2}, we obtain
$$
\rho(G)\leq\rho(G_1)\leq\rho(G_2)<\rho(G_*)=\rho(K_1\vee(K_{n-b-4}\cup K_3\cup bK_1)),
$$
which contradicts $\rho(G)\geq\rho(K_1\vee(K_{n-b-4}\cup K_3\cup bK_1))$. \hfill $\Box$

\section{The proof of Theorem 1.4}

\noindent{\it Proof of Theorem 1.4.} Suppose that a connected $\frac{1}{k-2}$-binding graph $G$ has no spanning $k$-tree. Using Lemma 2.2, we obtain $c(G-S)\geq(k-2)|S|+3$ for some nonempty subset $S$ of $V(G)$.
Write $|S|=s$ and $c(G-S)=q$. Then $q\geq(k-2)s+3$. We use $C_1,C_2,\ldots,C_q$ to denote the $q$ components in $G-S$ and write $|C_i|=n_i$ for $1\leq i\leq q$. Without loss of generality, let $n_q\geq n_{q-1}\geq\cdots\geq n_1$.
We are to verify $n_{(k-2)s+1}\geq2$. Otherwise, we have $n_i=1$ for $1\leq i\leq(k-2)s+1$. Set $X=V(C_1\cup C_2\cup\cdots\cup C_{(k-2)s+1})$. Obviously, $N_G(X)\subseteq S$ and $|X|=(k-2)s+1$. Thus, we infer
$$
\frac{|N_G(X)|}{|X|}\leq\frac{|S|}{|X|}=\frac{s}{(k-2)s+1}<\frac{1}{k-2},
$$
which is impossible because $G$ is $\frac{1}{k-2}$-binding. This leads to $n_i\geq2$ for $i\geq(k-2)s+1$. Clearly, $G$ is a spanning subgraph of $G_1=K_s\vee(K_{n_1}\cup\cdots\cup K_{n_{(k-2)s+1}}\cup K_{n_{(k-2)s+2}}\cup K_{n_{(k-2)s+3}})$,
where $n_{(k-2)s+3}\geq n_{(k-2)s+2}\geq n_{(k-2)s+1}\geq\cdots\geq n_1\geq1$, $n_{(k-2)s+1}\geq2$ and $\sum\limits_{i=1}^{(k-2)s+3}n_i=n-s$. According to Lemma 2.3, we obtain
\begin{align}\label{eq:4.1}
\rho(G)\leq\rho(G_1),
\end{align}
with equality holding if and only if $G=G_1$. Let $G_2=K_s\vee(K_{n-(k-1)s-4}\cup2K_2\cup(k-2)sK_1)$, where $n\geq(k-1)s+6$. By virtue of Lemma 2.6 (set $t=(k-2)s+3$ in Lemma 2.6), we get
\begin{align}\label{eq:4.2}
\rho(G_1)\leq\rho(G_2),
\end{align}
where the equality holds if and only if $G_1=G_2$. For $s=1$, we have $G_2=K_1\vee(K_{n-k-3}\cup2K_2\cup(k-2)K_1)$. By \eqref{eq:4.1} and \eqref{eq:4.2}, we conclude
$$
\rho(G)\leq\rho(K_1\vee(K_{n-k-3}\cup2K_2\cup(k-2)K_1)),
$$
with equality if and only if $G=K_1\vee(K_{n-k-3}\cup2K_2\cup(k-2)K_1)$, a contradiction. Next, we consider $s\geq2$.

Recall that $G_2=K_s\vee(K_{n-(k-1)s-4}\cup2K_2\cup(k-2)sK_1)$. Then
\begin{align*}
2e(G_2)=&2\binom{n-(k-2)s-4}{2}+8s+4+2(k-2)s^{2}\\
=&k(k-2)s^{2}+(9k-10-2(k-2)n)s+n^{2}-9n+24.
\end{align*}
Combining this with Lemma 2.4, we conclude
\begin{align}\label{eq:4.3}
\rho(G_2)\leq&\sqrt{2e(G_2)-n+1}\nonumber\\
=&\sqrt{k(k-2)s^{2}+(9k-10-2(k-2)n)s+n^{2}-10n+25}.
\end{align}
Let $\varphi(s)=k(k-2)s^{2}+(9k-10-2(k-2)n)s+n^{2}-10n+25$. Recall that $s\geq2$ and $n\geq(k-1)s+6$. Thus, we deduce $2\leq s\leq\frac{n-6}{k-1}$. By $k\geq3$ and $n\geq2k+22$, we possess
\begin{align*}
\varphi(2)-\varphi\left(\frac{n-6}{k-1}\right)=&\frac{(k-2)^{2}n^{2}-(4k^{3}-7k^{2}-11k+26)n+4k^{4}+2k^{3}-18k^{2}+8k+40}{(k-1)^{2}}\\
\geq&\frac{(k-2)^{2}(2k+22)^{2}-(4k^{3}-7k^{2}-11k+26)(2k+22)+4k^{4}+2k^{3}-18k^{2}+8k+40}{(k-1)^{2}}\\
=&\frac{306k^{2}-1386k+1404}{(k-1)^{2}}\\
\geq&\frac{306\times 9-1386\times 3+1404}{(k-1)^{2}}\\
=&0,
\end{align*}
which yields that
\begin{align}\label{eq:4.4}
\varphi(s)\leq\max\left\{\varphi(2),\varphi\left(\frac{n-6}{k-1}\right)\right\}\leq\varphi(2)
\end{align}
for $2\leq s\leq\frac{n-6}{k-1}$. From \eqref{eq:4.3}, \eqref{eq:4.4}, $k\geq3$ and $n\geq2k+22$, we deduce
\begin{align}\label{eq:4.5}
\rho(G_2)\leq&\sqrt{\varphi(s)}\leq\sqrt{\varphi(2)}\nonumber\\
=&\sqrt{n^{2}-(4k+2)n+4k^{2}+10k+5}\nonumber\\
=&\sqrt{(n-k-3)^{2}-((2k-4)n-3k^{2}-4k+4)}\nonumber\\
\leq&\sqrt{(n-k-3)^{2}-((2k-4)(2k+22)-3k^{2}-4k+4)}\nonumber\\
=&\sqrt{(n-k-3)^{2}-(k^{2}+32k-84)}\nonumber\\
<&n-k-3.
\end{align}
Notice that $K_{n-k-2}$ is a proper subgraph of $K_1\vee(K_{n-k-3}\cup2K_2\cup(k-2)K_1)$. According to Lemma 2.3, $\rho(K_1\vee(K_{n-k-3}\cup2K_2\cup(k-2)K_1))>\rho(K_{n-k-2})=n-k-3$. Together with
\eqref{eq:4.1}, \eqref{eq:4.2} and \eqref{eq:4.5}, we conclude
$$
\rho(G)\leq\rho(G_1)\leq\rho(G_2)<\rho(K_1\vee(K_{n-k-3}\cup2K_2\cup(k-2)K_1)),
$$
which contradicts $\rho(G)\geq\rho(K_1\vee(K_{n-k-3}\cup2K_2\cup(k-2)K_1))$. \hfill $\Box$

\section*{Declaration of competing interest}

The authors declare that they have no known competing financial interests or personal relationships that could have appeared to influence the work reported in this paper.

\section*{Data availability}

No data was used for the research described in the article.

\section*{Acknowledgments}

This work was supported by the Natural Science Foundation of Jiangsu Province (Grant No. BK20241949). Project ZR2023MA078 supported by Shandong Provincial Natural Science Foundation.

\end{document}